\documentstyle[11pt,twoside]{article}

\begin{document}
\overfullrule=0pt
\hsize=12,8cm
\vsize=19cm
\hoffset=1cm

\font\tenBbb=msbm10 \font\sevenBbb=msbm7 
\font\elevenBbb=msbm10 at 11pt
\font\twelveBbb=msbm10 at 12pt
\newfam\Bbbfam
\textfont\Bbbfam=\twelveBbb \scriptfont\Bbbfam=\sevenBbb
\def\Bbb{\fam\Bbbfam\tenBbb}

\font\teneufm=eufm10 \font\seveneufm=eufm7 
\font\twelveeufm=eufm10 at 12pt
\newfam\eufmfam
\textfont\eufmfam=\twelveeufm \scriptfont\eufmfam=\seveneufm
\def\frak#1{{\fam\eufmfam\relax#1}}
\let\goth\frak

\newcommand{\G}{\raisebox{0.1mm}{$\Gamma$}}

\newtheorem{De}{Definition}[section]
\newtheorem{Th}[De]{Theorem}
\newtheorem{Pro}[De]{Proposition}
\newtheorem{Le}[De]{Lemma}
\newtheorem{Ex}[De]{Example}
\newtheorem{Co}[De]{Corollary}

\newtheorem{dummy}[De]{}

\newcommand{\Def}[1]{\begin{De}#1\end{De}}
\newcommand{\Thm}[1]{\begin{Th}#1\end{Th}}
\newcommand{\Prop}[1]{\begin{Pro}#1\end{Pro}}
\newcommand{\Exa}[1]{\begin{Ex}#1\end{Ex}}
\newcommand{\Lem}[1]{\begin{Le}#1\end{Le}}
\newcommand{\Cor}[1]{\begin{Co}#1\end{Co}}

\def \Hom{\mathop{\rm Hom}\nolimits}

\def \Xi{\mathop{\rm {\cal F}({\sf as})}\nolimits}

\def\N{\noindent \par}
\def \n{\mathop{\underline n }\nolimits}
\def \m{\mathop{\underline m }\nolimits}
\def \k{\mathop{\underline k }\nolimits}
\def \g{\mathop{\Gamma }\nolimits}
\def\S{\smallskip \par}
\def\M{\medskip \par}

\def\B{\bigskip \par}
\def\BB{\bigskip \bigskip \par}

\def\t{\otimes }

\centerline{\bf { ON THE PROP CORRESPONDING TO BIALGEBRAS} }
\bigskip
\bigskip

\centerline{\bf {  By Teimuraz PIRASHVILI}}

\bigskip
 
\bigskip
\noindent {\sc R\'esum\'e}. Un PROP ${\bf A}$ est une cat\'egorie
mono\"{\i}dale  sym\'etrique stricte avec 
la propri\'et\'e suivante (cf \cite{Mac}): les 
objets de ${\bf A}$ sont les nombres naturels et l'op\'eration mono\"{\i}dale 
est l'addition sur les objets. Une alg\`ebre sur ${\bf A}$ est un foncteur
mono\"{\i}dal strict de ${\bf A}$  vers la cat\'egorie tensorielle
 ${\sf Vect}$ des espaces vectoriels sur un corps commutatif $k$. 
On construit le PROP $Q{\cal F}({\sf as})$ et on montre que
les alg\`ebres sur  $Q{\cal F}({\sf as})$  sont exactement les big\`ebres.

\bigskip

\section{Introduction} 
A PROP is a permutative  category $({\bf A }, 
\Box)$,  whose  set
of objects is the set of natural numbers and on 
objects the monoidal structure
is given by the addition.  An ${\bf A}$-{\it algebra} is a symmetric 
strict monoidal functor
to the tensor 
category of vector spaces.

It is well-known that there exists a PROP whose category of algebras is 
equivalent
to the category of bialgebras (= associative and 
coassociative  bialgebras). In \cite{Markl} there
is a description of this PROP in  terms of generators 
and relations. Here we give a more explicit construction of the same object. 
Our construction uses the Quillen's 
$Q$-construction for double categories given 
in \cite{FL}.

The paper is organized as follows: In Section 2 we recall the definition of 
PROP and show how to obtain commutative algebras as  ${\cal F}$-algebras.
Here  ${\cal F}$ is the PROP 
 of finite sets. In the next section following to \cite{FT} we construct the 
PROP of noncommutative sets denoted by ${\cal F}({\sf as})$ and we show that 
${\cal F}({\sf as})$-algebras are exactely associative algebras. 
The material of the Sections 2 and 3 are well 
known to experts. In Section 4 we
generalize  the notion of Mackey functor for double categories and in Section
5 we describe our hero ${\cal Q}$${\cal F}({\sf as})$, which is the PROP, 
with the property that ${\cal Q}$${\cal F}({\sf as})$-algebras 
are exactely bialgebras. By definition of PROP the category  
${\cal Q}$${\Xi}$ encodes
the natural transformations $H^{\t n}\to H^{\t m}$ and 
relations between them. Here $H$ runs over all bialgebras. As a sample
we give the following application. For any bialgebra $H$, any natural
number  $n\in \Bbb N$ and any permutation $\sigma \in {\frak S} _n$, we let
$$\Psi ^{(n, \sigma)}: H\to H$$
be the composition $\mu^n\circ \sigma _*\circ \Delta^n:H\to H$, where
$\Delta ^n: H\to H^{\t n}$ is the $(n-1)$-th iteration of the comultiplication
$\Delta : H\to H\t H$, 
$\sigma _* : H^{\t n}\to H^{\t n}$ is induced
by the permutation $\sigma$, that is
$$\sigma _*(x_1\t \cdots \t x_n)= x_{\sigma 1}\t \cdots \t x_{\sigma n}$$
 and $\mu ^n:  H^{\t n}\to H$ is 
the $(n-1)$-th iteration of the multiplication
$\mu : H\t H \to H$.  
Moreover let
$ \Phi: {\frak S} _n \times {\frak S} _m \to {\frak S} _{nm}$
be the map constructed in Proposition 5.3. Then it is a consequence of our
discussion in Section 5, that for any permutations
$\sigma \in {\frak S} _n$ and $\tau \in {\frak S} _m$
one has the equality
$$
\Psi ^{(n, \sigma)}\circ \Psi ^{(m, \tau)}= \Psi ^{(nm,\Phi(\sigma, \tau))}.
$$ 
 Let us note that if  $\sigma$ is the identity, then
$\Psi ^{(n, id)}$ is nothing but the Adams operation \cite{adams} and hence
our formula gives the rule for the composition of Adams operations.

\section{Preliminaries on PROP's}
 Recall that a {\it symmetric monoidal 
category} is a category ${\bf S}$ with a unit $0\in {\bf S}$ and a bifunctor
$$\Box: {\bf S}\times {\bf S}\to {\bf S}$$
together with natural isomorphisms 
$$a_{X,Y,Z}:X\Box (Y\Box Z)\to (X\Box Y)\Box Z,$$
$$ l_X: X\Box 0\to X, r_X:0\Box X \to X, c_{X,Y}:X\Box Y\to Y\Box X $$
satisfying some coherent conditions (see \cite{t}). If in 
addition $a_{X,Y,Z}$, $l_X$, $r_X$  are identity morphism then, ${\bf S}$ is called a 
{\it permutative category}. If  ${\bf S}$ and ${\bf S_1}$ are symmetric monoidal 
categories, then a functor $ M: {\bf S} \longrightarrow {\bf S_1}$
is a {\it  symmetric
monoidal
functor} if there exist  isomorphisms
$$ u_{X,Y}:M(X) \Box M(Y) \longrightarrow M(X \Box Y) $$
satisfying  the usual associativity and unit coherence conditions 
(see for example \cite{t}).                    
A  symmetric monoidal functor is called {\it strict} if $ u_{X,Y}$ is identity for all 
$X, Y\in {\bf S}$.  According  to \cite{Mac} a 
PROP is a permutative  category $({\bf A }, 
\Box)$,  with the following property: ${\bf A }$ has a  set
of objects equal to the set of natural numbers and on objects the bifunctor 
$\Box$ is given by $m\Box n=m+n$.  An ${\bf A }$-{\it algebra} is a symmetric 
strict monoidal functor
from ${\bf A }$ to the tensor 
category ${\sf Vect}$ of vector spaces over a field $k$.

\B
\noindent {\bf Examples}. 1) Let ${\cal F}$ be the category of finite sets.
For any $n\geq 0$,
we let $\n$ be the set $\lbrace 1,...,n\rbrace$. Hence $\underline 0$ is the
empty set. We 
assume that the objects of $\cal F$ are the sets $\n, \ n\geq 0$. The disjoint union makes the
category ${\cal F}$ a PROP. It is well-known that {\it the category of algebras
over  
${\cal F}$  is equivalent to the category
of commutative and associative algebras with unit}. Indeed, if $A$ is a commutative algebra, then the functor
${\cal L}_*(A):{\cal F}\to {\sf Vect}$ is a  $\cal F$-algebra. 
Here the functor ${\cal L}_*(A)$ is
given by
$${\cal L}_*(A)(\n)=  A^{\t n}.$$
For any  map $f:\n\to \m$, the action of $f$ on ${\cal L}_*(A)$ is given by
$$f_*(a_1\t \cdots \t a_n):\ = b_1\t \cdots \t b_m,  $$
where 
$$b_j=\prod _{f(i)=j}a_i, \ j=1,\cdots ,n. $$
Conversely, assume $T$ is a ${\cal F}$-algebra. We let $A$ be the value of $T$
on $\underline 1$.  The unique map $\underline 2\to \underline 1$ yields a homomorphism
$$\mu:A\t A\cong T(  \underline 2)\to T( \underline 1)= A.$$ 
On the other hand the unique map
 $\underline 0\to  \underline 1$ yields  a  
homomorphism $\eta:k=T(\underline 0)\to T(\underline 1)=A$.
The pair $(\mu, \eta)$ defines  on $A$ a structure of commutative 
and associative  algebra with unit. One can
use the fact that $T$ is a symmetric strict monoidal functor to prove that $T\cong  {\cal L}_*(A)$.

\M
\noindent 2) Let us note that  the opposite of a PROP is still a PROP 
with the same $\Box$. Hence the disjoint union yields also a structure of 
PROP on ${\cal F}^{op}$.
{\it The category of  ${\cal F}^{op}$-algebras is equivalent to the category
of cocommutative and coassociative coalgebras with counit}. For any such coalgebra $C$ we let 
${\cal L}^*(C):{\cal F}^{op}\to {\sf Vect}$ be the corresponding  ${\cal F}^{op}$-algebra. On objects 
we still have ${\cal L}^*(C)(\n)=  C^{\t n}.$

\M
\noindent 3) We let $\Omega$ be the subcategory of ${\cal F}$, which has the same objects
as ${\cal F}$, but morphisms are surjections. Clearly $\Omega$ is a 
subPROP of ${\cal F}$ and $\Omega$-algebras are (nonunital) 
commutative algebras.

\M
\noindent 4) We let ${\bf Mon}$ be the category of 
 finitely generated free monoids, which is a PROP 
with respect to coproduct. Similarly the category
${\bf Abmon}$ of finitely generated free abelian monoids, the category
${\bf Ab}$ of finitely generated free abelian groups and the category
${\bf Gr}$ of finitely generated free groups are  PROP's with respect to coproducts. For the category of algebras over these PROP's see Theorem 5.2 and Remark 1 at the end of the paper.

\M
\noindent 5) Any algebraic theory in the sense of Lawvere \cite{W} gives 
rise to a PROP. This generalizes the examples from 1) and 4).

\B
In the next section we give a noncommutative generalization of Examples 1)-3).

\section{Preliminaries on noncommutative sets} 
 In this section following to \cite{FT} we introduce the 
PROP ${\bf \Xi}$ with property that ${\bf \Xi}$-algebras
are associative algebras with unit. As a category ${\bf \Xi}$ is described in
\cite{FT}, p.191 under the name ``symmetric category''. It is also 
isomorphic to the category $\Delta S$ 
considered in \cite{CH}, \cite{FL}.  Objects of ${\bf \Xi}$
are finite sets. So $Ob({\cal F})= Ob (\Xi)$. A morphism from
$\n$ to $\m$ is a map $f:\n\to \m$ together with  a total ordering on 
$f^{-1}(j)$ for all $j\in \m$. By abuse of notation we will denote
morphisms in $\Xi$ by $f, g$ etc. Moreover sometimes we write $\mid f\mid$ for 
the underlying map 
of $f\in \Xi$. We will also say that $f$ is a noncommutative lifting of a  map $\mid f \mid$. 
In order to define the composition in $\Xi$ we recall 
the definition of ordered union of ordered sets. 
Assume $\Lambda $ is a totally ordered set and for each $\lambda \in \Lambda$ a totally 
ordered set $X_{\lambda}$ is given. 
Then $X= \coprod X_{\lambda}$  is 
the disjoint union  of the sets $X_{\lambda}$ which is ordered as follows. 
If $x\in X_{\lambda }$ and 
$y\in X_{\mu}$,  then $x\leq y$ in $X$ iff $\lambda < \mu$ or $\lambda =\mu$ and 
$x\leq y$ in $X_{\lambda}$. 

If $f\in \Hom_{\Xi}(\n,\m)$ and $g\in 
\Hom_{\Xi}(\m,\k)$, then the composite $gf$ is  $\mid g\mid  \mid f\mid$ as a map,
 while the total ordering in $(gf)^{-1}(i)$,
$i\in \k$ is given by the identification
$$(gf)^{-1}(i)=\coprod _{j\in g^{-1}(i)}f^{-1}(j).$$

Clearly one has the forgetful functor $\Xi \to \cal F$.
A morphism $f$ in $\Xi$ is called a {\it surjection} if the map $\mid f\mid$ is a surjection.
An {\it elementary surjection} is a surjection $f:\n\to \m$ for which $n-m\leq 1$. 

Since any injective map has the unique noncommutative lifting, we see that the disjoint
union, which defines the symmetric monoidal category 
structure in ${\cal F}$ has the unique lifting
in $\Xi$. Hence $\Xi $ is a PROP. 

We claim that {\it the category  $\Xi $-algebras
is equivalent to the category
of associative algebras with unit}. The only point here is the following.
Let us denote  by 
$\prod _{i\in I}^{<}x_{i}$ the product of the elements $x_{i}\in A$ 
where 
$I$ is a finite totally ordered set and the ordering 
in the product follows to the  ordering $I$. 
Here $A$ is an associative algebra.
Now  we have a  $\Xi $-algebra
${\cal X}_*(A):{\Xi }\to {\sf Vect}$. Here the functor ${\cal X}_*(A)$ is
given by the same rule as ${\cal L}_*(A)$ in the previous section, but
to take $\prod ^{<}$ in the definition of $b_j$. For example, if 
$f:\underline 4 \to \underline 3$ is given by $f(1)=f(2)=f(4)=3$, $f(3)=1$
and the total ordering in $f^{-1}(3)$ is $2<4<1$ then 
$f_*:A^{\t 4}\to A^{\t 3}$ is nothing but $a_1\t a_2\t a_3\t a_4\mapsto a_3\t 1\t a_2a_4a_1$. 

 We let $\Omega({\sf as})$ be the subcategory of ${\cal F}({\sf as})$, which has the same objects
as ${\cal F}({\sf as})$, but morphisms are surjections. Clearly $\Omega({\sf as})$ is a subPROP of 
${\cal F}({\sf as})$ and  $\Omega({\sf as})$-algebras 
are (nonunital) associative algebras.

Quite similarly, for any coassociative coalgebra $C$ with 
counit one has a ${\Xi }^{op}$-algebra
${\cal X}^*(C):{\Xi }^{op}\to {\sf Vect}$  with ${\cal X}^*(C)(\n)= C^{\t n}$ and   
{\it the category ${\Xi }^{op}$-algebras is equivalent to the category
of coassociative coalgebras with counit}.

In order to put bialgebras in the picture we need the language of Mackey functors.

\section{On double categories and Mackey functors} 
Let us recall that
 a  {\it double category} consists of objects, a set 
of horizontal morphisms, a set of vertical morphisms and a set of bimorphisms 
satisfying natural  conditions  \cite{ehre} (see also \cite{FL}). If ${\bf D}$ is a 
double category, we let ${\bf D}^h$ (resp. ${\bf D}^v$) be the category of
objects and horizontal  (resp. vertical)  morphisms of ${\bf D}$.

A  {\it Janus functor} $M$ from  a  double category ${\bf D}$ to  $\sf Vect$
 is the following data

\smallskip
\noindent i) a covariant  functor $M_*:{\bf D}^h \to {\sf Vect}$

\smallskip
\noindent ii) a contravariant  functor $M_*:({\bf D}^v)^{op} \to {\sf Vect}$

\smallskip
\noindent such that for each object $S\in {\bf D}$ one has 
$M_*(S)=M^*(S)=M(S)$.  A {\it Mackey functor}
 $M=(M_*, M^*)$ from  a  double category ${\bf D}$ to ${\sf Vect}$  is a Janus functor 
$M$ from a  double category ${\bf D}$ to ${\sf Vect}$
such that for each bimorphism in ${\bf D}$
$$\matrix{
 \ \ \ \ \ \ \ \ \  \ &U &\ \buildrel f_1 \over
\longrightarrow  \ &S\ \ \cr
\ &\ \ & \  \ & \ \ \cr
{\alpha }=\ \ \ \ \ \  \ & \downarrow & \phi _1 \ \ \ \ \ \ \ \ \ \ \ \ \phi &\downarrow \cr
\ & \ \ & \  \ & \ \ \cr
\ \ \ \ \ \ \ \ \  \ & T &\ \buildrel f \over
\longrightarrow \ & V \ \ \cr
}$$
the following equality holds:
$$M^*(\phi)M_*(f)=M_*({f_1})M^*({\phi _1}).$$

\noindent {\bf Examples} 1) Let ${\bf C}$ be  a category with pullbacks. 
Then one has  a double category 
whose objects are the same  as ${\bf C}$. Moreover
$Mor ^v= Mor ^h=Mor ({\bf C})$, while bimorphisms are pullback diagrams in 
${\bf C}$. In this
case the notion of Mackey functors corresponds to pre-Mackey functors from
\cite{dress}. 
By abuse of notation we will still denote this double category by ${\bf C}$. 
 In what follows $\cal F$ is equipped with this double category structure.

\M
\noindent 2) Now we consider a double category, whose objects
are still finite sets, but $Mor ^v= Mor ^h=Mor(\Xi)$, where $\Xi$ was introduced in Section 3. 
By definition a bimorphism is 
a diagram in $\Xi$ 

$$\matrix{
 \ \ \ \ \ \ \ \ \  \ &U &\ \buildrel f_1 \over
\longrightarrow  \ &S\ \ \cr
\ &\ \ & \  \ & \ \ \cr
{\alpha }=\ \ \ \ \ \  \ & \downarrow & \phi _1 \ \ \ \ \ \ \ \ \ \ \ \ \phi &\downarrow \cr
\ & \ \ & \  \ & \ \ \cr
\ \ \ \ \ \ \ \ \  \ & T &\ \buildrel f \over
\longrightarrow \ & V \ \ \cr
}$$
such that the following holds: 

\smallskip
\noindent i) the image $\mid \alpha \mid$ of 
$\alpha$ in ${\cal F}$ 
is a pullback diagram of sets, 

\smallskip
\noindent ii) for all $x\in T$ the 
induced map $f _*: \phi _1^{-1}(x)\to
\phi ^{-1}(fx)$ is an isomorphism of ordered sets 

\smallskip
\noindent  iii) for all $y\in S$ the induced map 
$\phi _*: f _1^{-1}(y)\to f ^{-1}(\phi _1y)$ is an isomorphism of
ordered sets.

Let us note that for a bimorhism $\alpha$ in $\Xi$ in general $f\circ \phi _1 
\not = \phi \circ f_1$. By abuse of notation we will denote 
this double category by $\Xi$. It is 
different from a double category considered 
in  \cite{FL}, which is also associated to the  category $\Xi$.

One observes that for any arrows $f:T\to V$, $\phi:S\to V$ 
in ${\cal F}({\sf as})$ there exists a bimorphism $\alpha$ which has
$f$ and  $\phi$ as edges and it is unique up to natural isomorphism. Indeed,
as a set we take $U$ to be the pullback and then we lift set maps $f_1$ and $\phi _1$ to the noncommutative world according to the properties ii) and iii). 
Clearly such lifting exists and it is unique.

  \B
\noindent 3) We can also consider the double category $\Xi_{1}$ 
whose
objects are still finite sets, vertical arrows are set maps, while  
horizontal ones are morphisms from $\Xi$.
The bimorphisms are diagrams similar to the diagrams in Example 2) but
such that $\phi$  and  $\phi _1$
are set maps, while  $f$ and $f_1$
are morphisms from $\Xi$. Furthermore the conditions i) and iii)  
from the previous example hold. We need also a double category
 $\Xi_2$ which is defined similarly, but now 
vertical arrows are morphisms 
from $\Xi$ and horizontal ones are set maps.

We have the following diagram of double categories, where arrows are forgetful functors
$$\matrix{
 \ \ \ \ \ \ \ \ \ \ \  \ \ \ \   \ & &  \Xi_{1} \ & &\ \ &\cr
\ &\ \ & \  \ & \ \ \cr
\ \ \ \ \ \ \  \ \ \ \ \ \  \ \ \  \ & \nearrow &  &\searrow &&\cr
\ & \ \ & \  \ & \ \ \cr
\ \ \  \ \ \ \Xi & &  & &{\cal F}. \ \ & &\ \ \ \  \ \ \ \  \ \ \ \ \  \ \  (4.0)\cr
\ &\ \ & \  \ & \ \ \cr
\ \ \ \ \ \  \ \ \    \  \ & \searrow &  &\nearrow &&\cr
\ \ \ \ \ \ \ \ \ \ \ \ \  \ \   \ & &  \Xi_{2} \ & \ \ &&\cr
}$$
Let ${\bf D}$ be  one of the double categories considered in (4.0).
A bimorphism $\alpha$ is called {\it elementary} if both $f$ and $\phi$ are elementary surjections.
The following Lemma for ${\bf D}=\cal F$ was proved in \cite{BDFP}. The proof in other 
cases is quite
similar and hence we omit it. 

\Lem{ Let ${\bf D}$ be one of the double categories considered in {\rm (4.0)}.
Then a 
Janus functor $M$ is a Mackey functor iff the following two conditions hold

i) for any injection $g:A\to B$ one has $M^*(g)M_*(g)=id_A$

ii) for  any elementary bimorphism $\alpha $ one has
 $$M^*(\phi)M_*(f)=M_*({f_1})M^*({\phi _1}).$$
}

\Thm{ Let $V$ be a vector space, which is equipped 
simultaneously with the structure of associative
algebra with unit and coassociative coalgebra with counit. 
Then $V$ is a bialgebra iff
$${\cal X}(V)=({\cal X}_*(V), {\cal X}^*(V)): \Xi \to {\sf Vect}$$
 is a Mackey functor.}

{\it Proof}. One observes that the condition 1) of the previous lemma 
always holds. On the other hand
the diagram
$$\matrix{
 \ \ \ \ \ \ \ \ \  \ &\underline 4 &\ \buildrel p \over
\longrightarrow  \ &\underline 2\ \ \cr
\ &\ \ & \  \ & \ \ \cr
\alpha =\ \ \ \ \ \  \ & \downarrow & q \ \ \ \ \ \ \ \ \ \ \ \ f &\downarrow \cr
\ & \ \ & \  \ & \ \ \cr
\ \ \ \ \ \ \ \ \  \ &\underline 2   &\ \buildrel f \over
\longrightarrow \ & \underline 1 \ \ \cr
}$$
is a bimorphism. Here $f^{-1}(1)=\{1<2\}$, $p^{-1}(1)=\{1<2\}$, $p^{-1}(2)=
\{3<4\}$, $q^{-1}(1)=\{1<3\}$ and
$q^{-1}(2)=\{2<4\}$. Clearly $f_*: V^{\t 2}\to V$ is the multiplication $\mu$ 
on $V$ and $f^*:V\to V^{\t 2}$ is the
comultiplication $\Delta$ on $V$, while $p_*=(\mu\t \mu)\circ \tau_{2,3}$ 
and $q^*=\tau_{2,3}\circ \Delta \t \Delta$,
where $\tau_{2,3}:V^{\t 4}\to V^{\t 4}$ permutes the 
second and the third coordinates.
Hence  $V$ is a bialgebra iff the condition ii)  of the previous lemma holds for $\alpha$. Since both
${\cal X}_*(V)$ and $ {\cal X}^*(V)$ send disjoint union to 
tensor product the result follows from Lemma 4.1.
 
\bigskip
\noindent {\bf Addendum}. For a cocommutative bialgebra $C$ the Mackey functor ${\cal X}(C)$
factors through the double category $\Xi _1$, for a commutative bialgebra $A$  
the Mackey functor ${\cal X}(A)$
factors through $\Xi _2$
and in the case of commutative and cocommutative bialgebra $H$
one has the Mackey functor
  ${\cal L}(H):{\cal F}\to {\sf Vect}.$ 

\section{The construction of ${\cal Q}$${\Xi}  $}
Let ${\bf D }$ be one of the double categories considered in Examples 1)-3). 
Clearly categories ${\bf D}^v$ and ${\bf D}^h$ have the same class of 
isomorphisms, which we call {\it isomorphisms of ${\bf D}$}. 
We let ${\cal Q}$${\bf  D}$ be the category whose objects are finite sets, 
while the morphisms from $T$ to $S$ are equivalence classes of diagrams:
$$\matrix{
 \ \ \ \ \ \ \ \ \  \ &U &\ \buildrel f \over
\longrightarrow  \ &S\ \ \cr
\ &\ \ & \  \ & \ \ \cr
\ \ \ \ \ \  \ & \downarrow & \phi \ \ \ \ \ \ \ \ \ \ \ \  & \cr
\ & \ \ & \  \ & \ \ \cr
\ \ \ \ \ \ \ \ \  \ & T &\ \ &  \ \ \cr
}$$
Here $f\in {\bf D}^h$ is a horizontal morphism and $\phi\in {\bf D}^v$ 
is a vertical morphism. For simplicity such data 
will  be denoted by $ T\buildrel \phi \over \hookleftarrow U \buildrel f \over \longrightarrow S$.
Two diagrams $ T\buildrel \phi \over
\hookleftarrow U \buildrel f \over
\longrightarrow S$ and $ T\buildrel \phi_1 \over
\hookleftarrow U_1 \buildrel f_1 \over \longrightarrow S$ 
are equivalent if there exists a commutative diagram
$$\matrix{ T &\ \buildrel \phi \over \hookleftarrow  \ &\ U\ & \buildrel f \over \longrightarrow  & S  \ \cr
\ \ \ & \  \ & \ \ \cr  \parallel & &h\downarrow & &\parallel\cr \ \ \ & \  \ & \ \ \cr T & \buildrel \phi _1 \over
\hookleftarrow  \    &\ U_1\ &\buildrel f_1 \over \longrightarrow  & S \  \cr }
$$
such that  $h$ is an isomorphism. The composition of  
$ T\buildrel \phi \over \hookleftarrow U \buildrel f \over \longrightarrow S$  
and $ S\buildrel \psi \over \hookleftarrow V \buildrel g \over \longrightarrow R$  
 in ${\cal Q}$${\bf D}$ is by definition 
$ T\buildrel \psi _1\phi \over  \hookleftarrow W \buildrel gf_1 \over \longrightarrow R$, where 
$$\matrix{
 W &\ \buildrel f_1 \over
\longrightarrow  \ &V\ \ \cr
\ \ \ & \  \ & \ \ \cr
\downarrow & \psi_1 \ \ \ \ \ \ \ \ \ \ \ \ \psi&\downarrow \cr
\ \ \ & \  \ & \ \ \cr
 U &\ \buildrel f \over
\longrightarrow \ & S. \ \ \cr
}$$
is a bimorphism in ${\bf  D}$.  One easily checks that ${\cal Q}$${\bf  D}$ is a category and for any object 
$S$ the  diagram $  S \buildrel 1_S \over \hookleftarrow S \buildrel 1_S \over 
\longrightarrow S$ is an identity morphism in ${\cal Q}$${\bf  D}$.

Clearly the disjoint union yields a  structure of 
PROP on ${\cal Q}$${\bf  D}$ and $\underline 0$ is  not only a unit object
 with respect to this monoidal structure, but also a zero object.   

For a horizontal morphism $f:S\to T$ in ${\bf D}$ we 
let $i_*(f):S\to T$ be the following morphism in ${\cal Q}$${\bf  D}$:
$$S\buildrel 1_S \over\hookleftarrow S\buildrel f \over \longrightarrow   T.$$ 
Similarly, for a vertical morphism $\phi :S\to T$  we let $i^*(f):T\to S$ be the 
following morphism in ${\cal Q}$${\bf  D} $:
$$T\buildrel f \over\hookleftarrow  S\buildrel 1_S \over \longrightarrow   S.$$  
In this way one obtains the morphisms of PROP's:
 $i_*:{\bf D }\to {\cal Q}$${\bf D}$ and 
$i^*:{\bf D}^{op}\to {\cal Q}$${\bf D}$. 

\B
{\bf Remark}. The construction of $\cal Q$${\bf D}$ is a particular case of 
the generalized Quillen $Q$-constuction \cite{Q} considered by
Fiedorowicz and   Loday in \cite{FL}. The following lemma is a variant 
of a result of \cite{L}.

\Lem{ The category of Mackey functors from  ${\bf D}$ to ${\sf Vect}$ 
is equivalent to the category of functors $M:\cal Q$${\bf D} \to {\sf Vect}$.
} 

\smallskip
{\it Proof}. Let $M:{\cal Q}$${\bf D}\to {\sf Vect}$ be a functor. 
For any arrow $f:S\to T$ we put $M_*(f):= M(i_*(f))$ and $M^*(f):= M(i^*(f))$. 
In this way we get a Mackey functor on ${\bf D}$. Conversely, if $M$ is  a Mackey 
functor on ${\bf D}$, then we put
$$M(S \buildrel g \over \hookleftarrow  V \buildrel f \over \longrightarrow T)=M_*(f)M^*(g).$$
One easily shows that in this way we get a covariant 
functor $\cal Q$${\bf D}$ to ${\sf Vect}$ and the proof is finished.

By applying the $Q$-construction to the diagram (4.0) one obtains the
following (noncommutative) diagram of PROP's:
$$\matrix{
 \ \ \ \ \ \ \ \ \ \ \  \ \ \ \   \ & &  {\cal Q}(\Xi_{1}) \  &\ \ &\cr
\ &\ \ & \  \ & \ \ \cr
\ \ \ \ \ \ \  \ \ \ \ \ \  \ \ \  \ & \nearrow &  &\searrow &\cr
\ & \ \ & \  \ & \ \ \cr
\ \ \  \ \ \ {\cal Q}(\Xi) & &  & &{\cal Q}({\cal F}) \ \ & \cr
\ &\ \ & \  \ & \ \ \cr
\ \ \ \ \ \  \ \ \    \  \ & \searrow &  &\nearrow &\cr
\ \ \ \ \ \ \ \ \ \ \ \ \  \ \   \ & & {\cal Q}( \Xi_{2}) \ & \ \ &\cr
}$$
The following theorem gives the identification of the terms involved in 
the diagram, except for ${\cal Q}(\Xi)$.

\Thm{ i) The category of ${\cal Q}({\Xi}) $-algebras 
is equivalent to the category of
bialgebras.

\smallskip

ii) The category ${\cal Q}({\Xi _1}) $-algebras is equivalent to the 
category of cocommutative
bialgebras and ${\cal Q}({\Xi _1}) $ is isomorphic to the PROP 
${\bf Mon}^{op}$.

\smallskip

iii) The category of ${\cal Q}({\Xi _2}) $-algebras is equivalent 
to the category of commutative
bialgebras and ${\cal Q}({\Xi _2}) $ is isomorphic to the PROP 
${\bf Mon}$.

\smallskip
iv) The category of ${\cal Q}({\cal F} )$-algebras is equivalent to the 
category of cocommutative
and commutative bialgebras and ${\cal Q}({\cal F} )$ is isomorphic to the PROP  ${\bf Abmon}$.
}

\S
{\it Proof}. Theorem 4.2  together with Lemma 5.1 shows that any bialgebra $V$ gives rise to ${\cal X}(V)$-algebra. Conversely assume $M$ is a 
${\cal Q}({\Xi}) $-algebra and let $V=M(\underline 1)$.
Then $M\circ i_*$ is a $\Xi$-algebra and $M\circ i^*$ is a 
$\Xi ^{op}$-algebra. Thus $M$ carries  natural
structures of associative algebra and coassociative coalgebra. 
Since $M=  (M\circ i_*, M\circ i^*)$ is
a Mackey functor on $\Xi$, it follows from Theorem 4.2
that $V$ is indeed a bialgebra. 
To prove the remaining parts of the theorem, let us observe that
$({\cal Q}({\Xi _2}) )^{op} \cong {\cal Q}({\Xi _1}) $, where equivalence is
identity on objects and sends 
$ T\buildrel \phi \over \hookleftarrow U \buildrel f \over \longrightarrow S$ to
$ S\buildrel f \over \hookleftarrow U \buildrel \phi \over \longrightarrow T$.
We now show that
 ${\cal Q}({\Xi _2}) \cong {\bf Mon}$.
 The main observation here is the fact that if $f:X\to S_1 \coprod S_2$ is a
 morphism in $\Xi$ then $f=f_1\coprod f_2$ in the category $\Xi$, where $f_i$ as a map is the 
restriction of $f$ on $f^{-1}(S_i), i=1,2$. 
Since $f_i^{-1}(y)=f^{-1}(y)$ for all $y\in f^{-1}(S_i)$ we can take the same
total ordering in $f_i^{-1}(y)$ to turn $f_i$ into a morphism in $\Xi$. 
A conclusion of this observation is the
fact that  disjoint union defines
 not only a symmetric monoidal category structure but it is the coproduct  in 
${\cal Q}({\Xi _2})$. 
Clearly  $\n$ is an $n$-fold coproduct of $\underline 1$.  
On the other hand, we may assume that the objects of ${\bf Mon}$ are 
natural numbers, while the set of morphisms from $k$ to $n$ is the same as 
$\Hom _{\bf Mon}(F_k, F_n)$, 
where $F_n$ is the free monoid on $n$ generators. This 
set can be identified with 
the set of $k$-tuples of words on
$n$ variables $x_1,\cdots , x_n$.  
Since ${\cal Q}({\Xi _2})$ and ${\bf Mon}$ are categories with finite
coproducts and any object in both categories is a  coproduct of some copies of $\underline 1$, we need 
only to identify the set of morphisms originating from $\underline 1$.
A morphism $\underline 1\to \n$ in 
${\cal Q}({\Xi _2})$ is a diagram 
$ \underline 1\buildrel 
\phi \over \hookleftarrow U \buildrel f \over \longrightarrow \n$,
where $\phi$ is a map of noncommutative sets.  We can associate to this
morphism a word $w$ of length $m$ on $n$ variables $x_1, \cdots , x_n$. 
Here $m= Card (U)$
 and the $i$-th place of $w$ is $x_{f(y_i)}$, where 
$U=\{ y_1< \cdots <y_m\}$. In this way
one sees immediately that this correspondence
defines the equivalence of categories ${\cal Q}({\cal F}({\sf as})_2)
\cong {\bf Mon}$. We refer the reader to \cite{BDFP}
for the fact that
${\cal Q}({\cal F}) $ is equivalent to ${\bf Abmon}$. Argument in this case
is even simpler than the previous one and  can be sketched as follows. Since
the PROP ${\cal Q}({\cal F}) $
is isomorphic to its opposite
disjoint union yields not only  the coproduct in   ${\cal Q}({\cal F}) $ but 
also  the product. Next, morphisms  $\underline
1 \to  \underline 1$ in ${\cal Q}({\cal F})$ are diagrams of maps
$\underline 1 \leftarrow U \to \underline 1$, whose equivalence class is
 completely determined by
the cardinality of $U$. This gives identification of morphisms from $\underline
1 \to \underline 1$ with natural numbers and the proof is done.

Thus the above diagram of PROP's is equivalent to the diagram 
$$\matrix{
 \ \ \ \ \ \ \ \ \ \ \  \ \ \ \   \ & &  {\bf Mon}^{op} \  &\ \ &\cr
\ &\ \ & \  \ & \ \ \cr
\ \ \ \ \ \ \  \ \ \ \ \ \  \ \ \  \ & \nearrow &  &\searrow &\cr
\ & \ \ & \  \ & \ \ \cr
\ \ \  \ \ \ {\cal Q}(\Xi) & &  & &{\bf Abmon} \ \ & \cr
\ &\ \ & \  \ & \ \ \cr
\ \ \ \ \ \  \ \ \    \  \ & \searrow &  &\nearrow &\cr
\ \ \ \ \ \ \ \ \ \ \ \ \  \ \   \ & &  {\bf Mon}\ & \ \ &\cr
}$$
Here ${\bf Mon}\to {\bf Abmon}$ is given by abelization functor.
Let us note that ${\cal Q}(\Xi)$ and ${\bf Abmon}$ are self dual PROP's,
and the arrows are surjection on morphisms. If one looks at endomorphisms
of $\underline 1$ we see that the endomorphism monoid 
$End_{\bf C}(\underline 1)$ for ${\bf C}= {\bf Mon}^{op}, {\bf Mon}, 
{\bf Abmon}$ is isomorphic to the multiplicative monoid of natural numbers. 
This corresponds to the fact that the operations $\Psi ^{(n, \sigma)}$ from 
the introduction for commutative or cocommutative bialgebras are independent
of $\sigma$ and  $\Psi ^n \circ \Psi ^m= \Psi ^{nm}$ \cite{adams}.

The following proposition describes the endomorphism monoid 
$End_{\bf C}(\underline 1)$ for ${\bf C}=Q({\cal F}({\sf as})).$

 Let $n\in \Bbb N$ be a natural number and let $\sigma\in {\frak S} _n$ 
be a permutation. Here $ {\frak S} _n$ is the group of permutations on $n$ 
letters. We let $[\sigma]$ be the morphism $\n\to \underline 1$ in $\Xi$
corresponding to the ordering $\sigma(1)< \sigma(2)<\cdots <\sigma(n)$.
For example $[id_n]$, or simply $[id]$ denotes 
the morphism $\n\to \underline 1$ in $\Xi$
corresponding to the ordering $1< 2<\cdots < n$.  Moreover we 
let $(n, \sigma):\underline 1\to \underline 1$ be the morphism 
 in 
${\cal Q}({\Xi })$ corresponding to the diagram 
$ \underline 1\buildrel [\sigma] \over \hookleftarrow \n \buildrel [id] 
\over \longrightarrow 1$.

\Prop{ The monoid of 
endomorphisms of $\underline 1 \in {\cal Q}({\Xi})$ is isomorphic to the
monoid of pairs $(n, \sigma)$, where $\sigma\in {\frak S} _n$ and 
$n\in \Bbb N$, with the following 
multiplication
$$(n, \sigma)\circ (m,\tau)=(nm, \Phi(\sigma, \tau)).$$
Here 
$$ \Phi: {\frak S} _n \times {\frak S} _m \to {\frak S} _{nm}$$
is a map, which is defined by
$$\Phi(\sigma, \tau)(x)=\tau (p+1)+m(\sigma (q)-1), \ \ 1\leq x\leq nm,$$
where $x=pn+q, \ 1\leq q\leq n$ and $0\leq p\leq m-1$.
}

{\it Proof}. A morphism $\underline 1\to \underline 1$ in ${\cal Q}({\Xi})$
is a diagram 
$ \underline 1\buildrel \phi \over \hookleftarrow U \buildrel f 
\over \longrightarrow \underline 1$, where 
$\phi$ and $f$  are morphisms of noncommutative sets. Hence $U$ has two 
total orderings corresponding to  $\phi$ and $f$. We will identify $U$ to $\n$,
via  ordering corresponding to $f$. Here $n$ is the cardinality of $U$. 
We denote the first (resp. the second, $\cdots $) element in the ordering corresponding to
$\phi$ by $\sigma (1)$ (resp. $\sigma (2), \cdots$ ). In this way we get 
a permutation $\sigma\in {\frak S} _n$. Thus any morphism 
$\underline 1\to \underline 1$ in ${\cal Q}({\Xi})$ is of the
form $(n, \sigma)$. In order to identify the composition law it 
is enough to note the following two facts: 

i) The  diagram 
$$\matrix{ {\underline {nm}} &\ \buildrel f \over \longrightarrow
  \ &\ \n\ &   \ \cr
\ \ \ & \  \ & \ \ \cr
g\downarrow \ \ & \  \ \  & \ \ \ \downarrow [\sigma]\cr \ \ \ & \  \ & \ \ \cr 
\m & \buildrel [id] \over \longrightarrow
 \     \ &\underline 1 \  \cr }
$$
is a bimorphism in ${\cal Q}({\Xi})$. Here $f$ and $g$ are given by
$$f^{-1}(j)=\{1+(j-1)m< 2+(j-1)m< \cdots < (m-1)+(j-1)m <jm\},$$
$$g^{-1}(i)=\{i+(\sigma (1)-1)m<i+(\sigma (2)-1)m
<\cdots <i+(\sigma (n)-1)m\},$$
for $i\in \m$ and $j\in \n$. 

\smallskip
ii) One has $[\Phi(\sigma, \tau)]=[\tau]\circ g$ and $[id_{\n}]\circ f=[id_{\underline {nm}}].$

We now give an alternative description of the function  $\Phi$. Let
$$\gamma : {\frak S}(n)\times  {\frak S}(m_1)\times \cdots\times {\frak S}(m_n)
\to {\frak S}(m_1+\cdots +m_n)\leqno (5.4)$$
be a map given by
$$\gamma(\sigma; \sigma_{m_1}, \cdots , \sigma_{m_n})= \sigma (m_1,\cdots ,m_n)\circ (\sigma_1 \coprod \cdots \coprod \sigma_{m_n}),$$
where $\sigma (m_1,\cdots ,m_n) $ permutes the $n$ blocks according to $\sigma$. Moreover, for any integers $n$  and $m$ we let 
$$I:{\underline {nm}}\to \n\times \m$$
be the bijection corresponding to the following ordering of the Cartesian product:
$$ (i,j)< (s,t) \ {\rm iff}\ i<s \ {\rm or} \ i=s \ {\rm  and }\ j<t.$$
Similarly, we let 
$$II:{\underline {nm}}\to \n\times \m$$
be the bijection corresponding to the following ordering of the Cartesian product:
$$ (i,j)< (s,t) \ {\rm iff}\ j<t \ {\rm or} \ j=t \ {\rm  and }\ i<s.$$
Then we put $\Phi(n,m):= I^{-1}\circ II\in  {\frak S}_{nm}$. 
It is not too difficult to see that  $\Phi (n,m)=\Phi (1_{\n},1_{\m}) $ and
$$ \Phi(\sigma, \tau)= \Phi (n,m)\circ \gamma(\tau, \sigma , \cdots, \sigma).$$

\noindent {\bf Remarks}. 
1) It is well known that
the PROP corresponding to
cocommutative Hopf algebras is ${\bf Gr}^{op}$ (see next remark), 
the PROP corresponding to
commutative Hopf algebras is ${\bf Gr}$, while 
the PROP corresponding to commutative and
cocommutative Hopf algebras is ${\bf Ab}$. 
Of course the category of Hopf algebras are
also algebras over some PROP, which can be easily described via
generators and relations  \cite{Markl}. An explicit  
description of this particular   PROP will be the subject of a forthcoming paper. 

\S
\noindent 2) Let $A$ be a cocommutative Hopf algebra. Since $\t$ is a 
product in the category 
${\bf Coalg}$ of cocommutative coalgebras, $A$ is a group
object in this category. On the other hand any group object in any category ${\bf A}$ with finite products
gives rise to the model in ${\bf A}$ 
of the algebraic theory of groups in the sense of Lawvere \cite{W}. But the 
 algebraic theory of groups is nothing but  ${\bf Gr}^{op}$ and hence we have the functor 
${\cal X}(A):{\bf Gr}^{op}\to {\bf Coalg}$, which
 assigns  $A^{\t n}$ to $<n>$. Here 
$<n>$ is a free group on $x_1, \cdots , x_n$. Moreover 
it assigns $\mu$ to the morphism $<1>\ \to \ <2>$ given by 
$ x_1\mapsto x_1x_2$. 
Similarly ${\cal X}(A)$ assigns $\Delta$ to
the homomorphism $<2>\ \to <1>\ $ given by 
$x_1, x_2\mapsto x_1$. 
Of course it assigns the antipode $S:A\to A$ to  
$x_1\mapsto x_1^{-1}$. Having these
facts in mind one easily describes the action of ${\cal X}(A)$ on more complicate morphisms. For example one
checks that ${\cal X}(A)$ assigns
$$(\mu,\mu)\circ (\mu,id,\mu,id)\circ (S,id_{A^{\t 4}})\circ \tau_{2,3}\circ 
(id_{A^{\t 3}},\Delta, id)\circ (\Delta ,\Delta ,id)$$
to the morphism $<2>\ \to \ <3>$ corresponding to the pair of words 
$(x_1^{-1}x_2x_1, x_1^2x_3)$. Here $\tau_{2,3}$ permutes the second and third
coordinates. Conversely any linear map $A^{\t n}\to A^{\t m}$ constructed 
using the structural data of a cocommutative Hopf algebra $A$ is comming in
this way. Hence to check whether a complicated diagram involving such
maps commutes it is enough to look to the corresponding diagram in ${\bf Gr}$,
which is usually simpler to handle.

\S
\noindent 3) It is well known that the morphism $\n\to \m$ in ${\bf Abmon}$
can be identified with $(m\times n)$-matrices over natural numbers. Under this 
identification the equivalence ${\cal Q}({\cal F})\cong {\bf Abmon}$ is given by assigning
the matrix whose $(i,j)$-component is the cardinality of $f^{-1}(j)\bigcap g^{-1}(i)$, $1\leq i\leq m, \ 1\leq j\leq n$ to the diagram
$ \underline n\buildrel f \over \hookleftarrow X \buildrel g 
\over \longrightarrow m$. It is less known that the morphisms $\n\to \m$ in
${\bf Mon}$
can be described via shuffles. In order to explain this connection let us 
start with particular case. 
Consider a word $x^2yxy^3x^2$ of 
bidegree $(5,4)$. It defines a morphism $ \underline 1\to \underline 2$ in  
${\bf Mon}$. One associates a $(5,4)$-shuffle $(1,2,4,8,9,3,5,6,7)$
to this word, whose first five values are just the numbers of places where $x$
lies. Similarly morphisms $\n \to \m$ in ${\bf Mon}$ are in 1-1-correspondence with collections  $\{ A=(a_{ij}), (\varphi _1, \cdots, \varphi _n)\}$, where
$A$ is an $(m\times n)$-matrix  over natural numbers and $ \varphi _i$ is 
a $(a_{i1}, \cdots , a_{im})$-shuffle, $i=1,\cdots , n$. The functor ${\bf Mon \to Abmon}$
corresponds to forgetting the shuffles. Now combine this observation with
Proposition 5.3 to get the description of morphisms  $\n\to \m$ in 
${\cal Q}({\Xi})$ as collections  
$\{ A=(\alpha _{ij}), (\varphi _1, \cdots, \varphi _n)\}$, where
$\alpha =(a_{ij}, \sigma _{ij})$ and
 $a_{ij}$ is a natural number, while  $\sigma _{ij} \in {\frak S}_{a_{ij}}$ is a permutation and finally
$ \varphi _i$ is a $(a_{i1}, \cdots , a_{im})$-shuffle.

\S
\noindent 4) Recently Sarah Whitehouse (\cite{w}, \cite{cw}) 
defined the action of ${\frak S}_{k+1}$ on $A^{\t k}$ for 
any commutative or cocommutative Hopf algebra $A$. 
Actually she implicitly constructed the
group homomorphism 
$$\xi_k:{\frak S}_{k+1} \to {\frak G}_k,$$
 where ${\frak G}_k$ is the automorphism group of $<k>$. Then the action 
of $x\in {\frak S}_{k+1}$ on $A^{\t k}$ is obtained by 
applying the functor ${\cal X}(A)$ to $\xi _k(x)$. The
homomorphism $\xi_k$ is given by 
$$\sigma _1(x_1)=x_1^{-1}, \ \sigma_ 1(x_2)=x_1x_2, \ \sigma _1(x_i)=x_i, \ i\geq 2$$
$$\sigma _i(x_{i-1})=x_{i-1}x_i, \ \sigma _i(x_i)=x_i^{-1}, \ 
\sigma _i(x_{i+1})=x_ix_{i+1}, \ 
\sigma _i(x_j)=x_j,$$
for  $1< i<k, \ j\not =i-1,i,i+1$ and
$$\sigma _k(x_{k-1}) =x_{k-1}x_k,\  \sigma _k (x_k)=x_k^{-1}, \ \sigma _k(x_j)=x_j 
\ {\rm if }\ j<n-1.$$
Here $\sigma _i\in {\frak S}_{k+1}$ is the transposition 
$(i,i+1), \ 1\leq i\leq k$.
The homomorphisms $\xi _k$, $k\geq 0$ are restrictions of a functor 
$\xi:{\cal F}\to {\bf Gr}$, which is given as follows. For a set $X$ 
the group $\xi(X)$ is generated by symbols $<x,y>$, 
$x,y\in X$ modulo the realtions
$$<x,y><y,z>=<x,z>, \ x,y,z\in X.$$

\section{Generalization for operads}

Let ${\cal P}$ be an operad of sets \cite{oper}. Let us recall that
then $\cal P$ is a collection
of ${\frak S}_n$-sets ${\cal P}(n)$, $n\geq 0$ together with the composition
law 
$$\gamma : {\cal P}(n)\times  {\cal P}(m_1)\times \cdots\times {\cal P}(m_n)
\to {\cal P}(m_1+\cdots +m_n)$$
and an element $e\in {\cal P}(1)$ satisfying some associativity
and unite conditions \cite{oper}. We will assume that ${\cal P}(0)=*$.
Any set $X$ gives rise to an operad ${\cal E}_X$, for which  ${\cal E}_X(n)=
{\sf Maps}(X^n,X)$. A $\cal P$-algebra is a set $X$ together with a 
morphism of operads ${\cal P}\to {\cal E}_X$. We let $\cal P$-${\sf Alg}$ be 
the category of $\cal P$-algebras. The forgetfull functor 
$\cal P$-${\sf Alg}\to {\sf Sets}$ has the left adjoint functor
$F_{\cal P}:{\sf Sets}\to \cal P$-${\sf Alg}$ which is given by
$$F_{\cal P}(X)=\coprod _{n\geq 0}  {\cal P}(n)\times_{\frak S_n} X^n.$$
 We let ${\sf Free}({\cal P})$ be the full subcategory of 
$\cal P$-${\sf Alg}$ whose objects are $F_{\cal P}(\underline n)$, $n\geq 0$.

Now we introduce the category ${\cal F}({\cal P})$.
For any map $f:\underline n\to \underline m$ one puts 
$${\cal P}_f= \prod_{i=1}^m {\cal P}(\mid f^{-1}(i)\mid).$$
Here $\mid S\mid$ denotes the cardinality of a set $S$.
The category  ${\cal F}({\cal P})$ has the same objects as ${\cal F}$, while 
the morphisms from $\underline n\to \underline m$
in ${\cal F}({\cal P})$ are pairs $(f, \omega ^f)$, where $f:\underline n\to \underline m$
is a map and $\omega ^f= (\omega_1 ^f , \cdots ,\omega_m ^f)
\in {\cal P}_f$. If $(f, \omega ^f)$ and $(g, \omega ^g):\underline m\to \underline k$ 
are morphisms in ${\cal F}({\cal P})$ then the composition 
$(g, \omega ^g)\circ (f, \omega ^f)$ is a pair $(h, \omega ^h)$, where $h=gf$ 
and for each $1\leq i\leq k$ one has 
$$\omega_i^h=\gamma (\omega _i^g; \omega^f_{j_1},\cdots ,\omega^f_{j_s}).$$
Here $g^{-1}(i)=\lbrace j_1,\cdots , j_s\rbrace$.  This construction goes back
to May and Thomason \cite{MT}. 

One observes that if $\cal P={\sf as}$, then
${\cal F}({\cal P})$
is nothing but ${\cal F}({\sf as})$, while ${\sf Free}({\cal P})$ is equivalent
to the category of finitely generated free monoids. Here $\sf as$
is the operad given by ${\sf as}(n)={\frak S}_n$
for all $n\geq 0$ and $\gamma$ is the same as in (5.4). Thus ${\sf as}$-algebras 
are associative monoids. We now show how to generalize Theorem 5.2 ii) 
for arbitrary operads.

Let ${\cal F}({\cal P})_2$ be the double category, 
whose objects are sets, horisontal
arrows are set maps and vertical arrows are 
morphisms from ${\cal F}({\cal P})$.
Double morphisms are pulback diagrams of sets
$$\matrix{
 \ \ \ \ \ \ \ \ \  \ &U &\ \buildrel p \over
\longrightarrow  \ &S\ \ \cr
\ &\ \ & \  \ & \ \ \cr
\alpha =\ \ \ \ \ \  \ & \downarrow & g \ \ \ \ \ \ \ \ \ \ \ \ f &\downarrow \cr
\ & \ \ & \  \ & \ \ \cr
\ \ \ \ \ \ \ \ \  \ &T  &\ \buildrel q \over
\longrightarrow \ & U \ \ \cr
}$$
together with lifting of $g$ and $f$ in ${\cal F}({\cal P})$. Hence the 
elements $\omega ^f\in {\cal P}_f$ and $\omega ^g\in {\cal P}_g$ are given.
One requires that these elements are compatible  
$$\omega ^g_t=\omega ^f_{qt}, \ \ t\in T.$$ 
We claim that the category
${\cal Q}({\cal F}({\cal P})_2)$ and ${\sf Free }({\cal P})$ are equivalent.
On objects one assignes $F_{\cal P}({\underline n}) $ to $\underline n$.  
Both categories in the question posses finite coproducts and thus one needs 
only to identify
morphisms from $\underline 1$. Let
$ \underline 1\buildrel 
\omega \over \hookleftarrow m \buildrel f \over \to X$
be a morphism in ${\cal Q}({\cal F}({\cal P})_2)$. By definition $\omega \in {\cal P}(m)$
and $f\in X^n$. Thus it gives an element in $F_{\cal P}(X)$ and therefore a
morphism $F_{\cal P}({\underline 1})\to F_{\cal P}(X)$ in ${\sf Free}(X)$. It is clear
that in this way one obtains expected equivalence of categories.

Any set operad ${\cal P}$ gives rise to the linear operad $k[{\cal P}]$, which is
spanned on ${\cal P}$.
Clearly the disjoint 
union yields a structure of PROP
on ${\cal F}({\cal P})$ and ${\cal F}({\cal P})$-algebras
are nothing but $k[{\cal P}]$-algebras in 
the tensor category ${\sf Vect}$.

We leave as an exercise to the interested 
readers to show that the 
$Q$-construction of Section 5 and the notion of the bialgebra have the
canonical generalizations  for any operad ${\cal P}$.

\bigskip

\centerline{\bf Acknowledgments}

\bigskip

\noindent This work was written during my visit at the 
Sonderforschungsbereich der
Universit\"at Bielefeld. I would like to thank Friedhelm Waldhausen for
the invitation to  Bielefeld.  It is a pleasure to acknowledge various helpful discussions I had 
with V. Franjou and J.-L. Loday on the subject and also for invitations in 
Nantes and Le Pouliguen, where this work was started. Thanks to Ross Street 
on his comments.
The author  was partially supported by the grant 
INTAS-99-00817 and by the TMR network K-theory and algebraic groups, ERB FMRX CT-97-0107. 

\noindent 
A.M. Razmadze Math. Inst., Aleksidze str. 1,
\par \noindent
Tbilisi, 380093. Republic of Georgia \par
\noindent
{ email: pira@rmi.acnet.ge} 
\end{document}